\documentclass[11pt,reqno]{amsart}
\usepackage{amscd,amssymb,amsmath,amsthm}
\usepackage[arrow,matrix]{xy}
\usepackage{graphicx}
\usepackage{cite}
\usepackage{geometry}
\usepackage[colorlinks,urlcolor=blue]{hyperref}
\newtheorem{thm}[subsection]{Theorem}
\newtheorem{lemma}[subsection]{Lemma}
\newtheorem{pro}[subsection]{Proposition}

\numberwithin{equation}{section} 

\newcommand{\s}{{\sigma}}
\newcommand{\de}{{\xi}}

\def \s {\sigma}

\newcommand{\bea}{\begin{eqnarray}}
\newcommand{\eea}{\end{eqnarray}}

\begin{document}
\title[On positive fixed points of operator of Hammerstein type]{On positive fixed points of operator of Hammerstein type with degenerate kernel and Gibbs Measures}

\author{I. M. Mavlonov,  Kh. N. Khushvaktov, G.P. Arzikulov, F.H.Haydarov}

\address{I. M. Mavlonov \\ National University of Uzbekistan,
Tashkent, Uzbekistan.} \email {mavlonovismoil16@gmail.com}

\address{Kh. N. Khushvaktov \\ National University of Uzbekistan,
Tashkent, Uzbekistan.} \email{nuriddinh@gmail.com}

\address{G.P. Arzikulov \\ Tashkent State Technical University,
Tashkent, Uzbekistan.} \email{arzikulov79 @mai.ru}

\address{F.H.Haydarov \\ New Uzbekistan University, 54, Mustaqillik Ave., Tashkent, 100007, Uzbekistan.}

\address{AKFA University, 264, Milliy Bog street,  Yangiobod QFY, Barkamol MFY,
 Kibray district, 111221, Tashkent region, Uzbekistan.} \email{haydarov\_imc@mail.ru}

\begin{abstract} From \cite{re} it is known that ``translation-invariant Gibbs measures" of the model with an uncountable set of spin values can be described by positive fixed points of a nonlinear integral operator of Hammerstein type. In \cite{enh2015, MSSX} there are main results on positive fixed points of the operator of Hammerstein type with degenerate kernels, but it was not solved the existence of Gibbs measures corresponding to the founded fixed points for constructed kernels.
This paper is an investigation of the papers \cite{enh2015} and \cite{MSSX}. In this paper we construct new degenerate kernels of the Hammerstein operator by taking into account problems in the theory of Gibbs measure, i.e. each positive fixed point of the operator gives translational-invariant Gibbs measure.
\end{abstract}
\maketitle



\section{Introduction} \label{sec:intro}

Hammerstein equation covers a large variety of areas and is of much interest
to a wide audience due to the fact that it has applications in numerous
areas. Several problems that arise in differential equations (ordinary and partial),
for instance, elliptic boundary value problems whose linear parts possess
Green's function can be transformed into the Hammerstein integral equations.
Equations of the Hammerstein type play a crucial role in the theory of optimal
control systems and in automation and network theory (see e.g., Dolezale \cite{dolezale}).

There are some works devoted to fixed points of Hammerstein operator on cones. We can find main results on the existence and multiplicity of fixed points of Hammerstein equations (e.g., \cite{appl1, Cabada24,  margot}). On the other hand, we need to find new results on the uniqueness of fixed points of Hammerstein equations in cones. For instance, during last years, an increasing attention was given to
models with a \emph{uncountable} many spin values on a Cayley tree. In \cite{enh 2015}, \cite{re} Hamiltonian with a {\it
uncountable} set of spin values (with the set $[0, 1]$ of spin
values) on a Cayley tree $\Gamma^k$ was considered  and it was
showed that, the existence translation-invariant splitting Gibbs
measure of the Hamiltonian is equivalent to the existence of a
positive fixed point of Hammerstein type nonlinear integral
operator.   In \cite{re} for $k = 1$ (when the Cayley tree becomes a one-dimensional
 lattice $\mathbb{Z}$) it is shown that the integral equation has a unique solution, implying that there
 is a unique Gibbs measure.  For general $k\geq 2$,
  a sufficient condition is found under which a periodic Gibbs measure is unique (see \cite{ehr2013}). On the other hand,
  on the Cayley trees $\Gamma_{k}$ of order $k\geq 2$, the existence of phase transitions has been proven, see
   \cite{Ro, ehr2012}. We note that all of these papers are devoted to models with
    nearest-neighbor interactions. Also, in \cite{h, h1, 26} the splitting Gibbs
measures for four competing interactions (external field, nearest neighbor, second neighbors and triples of neighbors)
of models on $\Gamma_{2}$ are described and the fact that periodic Gibbs measure for the Hamiltonians with four
competing interactions are either \emph{translation-invariant} or\emph{ periodic with period two} is shown.
But, theorems on fixed points of Hammerstein operator on cone are impossible to use directly from known results devoted to problems of existence and uniqueness of fixed points of Hammerstein integral operator on cones. In present paper we construct new degenerate kernels of Hammerstein operator by taking into account problems in the theory of Gibbs measure.

\section{Preliminaries}

The Cayley tree $\Im^{k}=(V, L)$ of order $k \geq 1$ is an infinite tree, i.e. graph without cycles, each vertex of which has exactly $k+1$ edges. Here $V$ is the set of vertices of $\Im^{k}$ va $L$ is the set of its edges.


 Consider models where the spin takes values in the set $[0,1]$,
and is assigned to the vertices of the tree. For $A\subset V$ a
configuration $\sigma_A$ on $A$ is an arbitrary function $\sigma_A:A\to
\Phi$. Let $\Omega_A=[0,1]^A$ be the set of all configurations
on $A$. A configuration $\sigma$ on $V$ is defined as a function
$x\in V\mapsto\sigma (x)\in [0,1]$; the set of all configurations
is $\Omega:=[0,1]^V$. We consider all elements of $V$ are numerated (in any order) by the numbers: $0,1,2,3,...$. Namely, we can write $V=\{x_0, x_1, x_2, ....\}$.

 Let $\mathcal{X}_A$ be the indicator function. $\Omega$ can be considered as a metric space with respect to the metric $\rho: \Omega \times \Omega \rightarrow \mathbb{R}^{+}$ given by
$$
\rho\left(\left\{\sigma(x_n)\right\}_{x_n \in V},\left\{\sigma^{\prime}(x_n)\right\}_{x_n \in V}\right)=\sum_{n \geq 0} 2^{-n}\mathcal{X}_{\sigma(x_n)\neq\sigma^{\prime}(x_n)}
$$
(or any equivalent metric the reader might prefer, this metric taken from \cite{2}), and let $\mathcal{B}$ be the $\sigma$-field of Borel subsets of $\Omega$.

For each $m \geq 0$ let $\pi_m: \Omega \rightarrow[0,1]^{m+1}$ be given by $\pi_m\left(\sigma_0, \sigma_1, \sigma_2, ...\right)=\left(\sigma_0, \ldots, \sigma_m\right)$ and let $\mathcal{C}_m=\pi_m^{-1}\left(\mathcal{P}\left([0,1]^{m+1}\right)\right)$, where $\sigma_i:=\sigma(x_i)$ and $\mathcal{P}\left([0,1]^{m+1}\right)$ is the family of all subsets of $[0,1]^{m+1}$ (Cartesian product of $[0,1]$). Then $\mathcal{C}_m$ is a field and each of the sets in $\mathcal{C}_m$ is open and closed set in the metric space $(\Omega, \rho)$; also $\mathcal{C}_m \subset \mathcal{C}_{m+1}$. Let $\mathcal{C}=\bigcup_{m \geq 0} \mathcal{C}_m$; then $\mathcal{C}$ is a field (the field of \textbf{cylinder sets}) and each of the sets in $\mathcal{C}$ is both open and closed.
Denote  $\mathcal{S}(\mathcal{C})$ - the smallest sigma field containing $\mathcal{C}$. Every element of $\mathcal{S}(\mathcal{C})$ is called ``\textbf{measurable cylinder}".

Let us consider a formal Hamiltonian:
\begin{equation}\label{e1.1}
 H(\sigma)=-J\sum_{\langle x,y\rangle\in L}
\de_{\sigma(x), \sigma(y)}, \,\ \sigma\in\Omega_{V}
\end{equation}
where $J \in R\setminus \{0\}$ and $\de: (u,v)\in [0,1]^2\to
\de_{uv}\in R$ is a given bounded, measurable function. As
usual, $\langle x,y\rangle$ stands for the nearest neighbor
vertices.

Let $h:\;x\in V\mapsto h_x=(h_{t,x}, t\in [0,1]) \in R^{[0,1]}$ be
mapping of $x\in V\setminus \{x^0\}$.  Given $n=1,2,\ldots$,
consider the probability distribution $\mu^{(n)}$ on
$\Omega_{V_n}$ defined by
\begin{equation}\label{e2}
\mu^{(n)}(\sigma_n)=Z_n^{-1}\exp\left(-\beta H(\sigma_n)
+\sum_{x\in W_n}h_{\sigma_n(x),x}\right),
\end{equation}
 Here, as before, $\sigma_n:x\in V_n\mapsto
\sigma(x)$ and $Z_n$ is the corresponding partition function:
\begin{equation}\label{e3} Z_n=\int_{\Omega_{V_n}}
\exp\left(-\beta H({\widetilde\sigma}_n) +\sum_{x\in
W_n}h_{{\widetilde\sigma_n}(x),x}\right)
\lambda_{V_n}({d\widetilde\s_n}).
\end{equation}

Let $\Lambda\in \mathcal{N}$ and $\Delta\subset \Lambda$. If $\mu_{\Lambda}$ is a  measure on $\mathcal{B}_{\Lambda}$, the projection of $\mu_{\Lambda}$ on $\mathcal{B}_{\Delta}$ is measure $\pi_{\Delta}\left(\mu_{\Lambda}\right)$ on $\mathcal{B}_{\Delta}$ defined by
$$
\left[\pi_{\Delta}\left(\mu_{\Lambda}\right)\right](B)=\mu_{\Lambda}\left\{\sigma \in \Omega_{\Lambda}: \ \sigma |_{\Delta} \in B\right\}, \ B \in \mathcal{B}_{\Delta}.
$$
Similarly, if $\mu$ is a  measure on $\mathcal{B}$, the projection of $\mu$ on $\mathcal{B}_{\Lambda}$ is defined by
$$
\left[\pi_{\Lambda}(\mu)\right](B)=\mu\left\{\sigma \in \Omega: \sigma_{\Lambda} \in B\right\}=\mu(\bar{\sigma} |_{\Lambda}=\sigma_{\Lambda} : \sigma_{\Lambda}\in B), \quad B \in \mathcal{B}_{\Lambda}.$$
The following theorem is known:
\begin{thm}\label{Kol}\cite{1}\textbf{(Kolmogorov Extension Theorem)} For each $t$ in the arbitrary index set $T$, let $\Omega_t$ be a complete, separable metric space, with $\mathcal{F}_t$ the class of Borel sets (the $\sigma$-field generated by the open sets).

Assume that for each finite nonempty subset $v$ of $T$, we are given a probability measure $P_v$ on $\mathcal{F}_v$. Assume the $P_v$ are consistent, that is, $\pi_u\left(P_v\right)=P_u$ for each nonempty $u \subset v$.

Then there is a unique probability measure $P$ on $\mathcal{F}=\prod_{t \in T} \mathcal{F}_t$ such that $\pi_v(P)=P_v$ for all $v$.
\end{thm}

The probability distributions $\mu^{(n)}$ are compatible if for
any $n\geq 1$ and $\sigma_{n-1}\in\Omega_{V_{n-1}}$:
\begin{equation}\label{e4}
\pi_{V_{n-1}}\left(\mu^{(n)}\right)=\mu^{(n-1)}
\end{equation}
Then by Kolmogorov extension theorem, there exists a unique
measure $\mu$ on $\Omega_V$ such that, for any $n$ and
$\sigma_n\in\Omega_{V_n}$, $\mu \left(\left\{\sigma
\Big|_{V_n}=\sigma_n\right\}\right)=\mu^{(n)}(\sigma_n)$.

The measure $\mu$ is called {\it splitting Gibbs measure}
corresponding to Hamiltonian (\ref{e1.1}) and function $x\mapsto
h_x$, $x\neq x^0$.

\begin{pro}\label{p1}\cite{re} {\it The probability distributions
$\mu^{(n)}(\sigma_n)$, $n=1,2,\ldots$, in} (\ref{e2}) {\sl are
compatible iff for any $x\in V\setminus\{x^0\}$ the following
equation holds:
\begin{equation}\label{e5}
 f(t,x)=\prod_{y\in S(x)}{\int_0^1\exp(J\beta\de_{tu})f(u,y)du \over \int_0^1\exp(J\beta{\de_{0u}})f(u,y)du}.
 \end{equation}
Here, and below  $f(t,x)=\exp(h_{t,x}-h_{0,x}), \ t\in [0,1]$ and
$du=\lambda(du)$ is the Lebesgue measure.}
\end{pro}

Note, that the analysis of solutions to (\ref{e5}) is not easy.
It's difficult to give a full description for the given potential
function $\de_{t,u}$.

Let $\xi_{tu}$ is a continuous function. We put
$$C^+[0,1]=\{f\in C[0,1]: f(x)\geq 0\}, \,\ C_0^+[0,1]=C^+[0,1]\setminus
\{\theta\equiv 0\}.$$

Define the operator $R_{k}:C^{+}_{0}[0,1]\rightarrow
C^{+}_{0}[0,1]$ by
$$(R_{k}f)(t)=\left({\int_0^1K(t,u)f(u)du\over \int_0^1
K(0,u)f(u)du}\right)^k, \,\
k\in\mathbb{N},$$\\
where $K(t,u)=\exp(J\beta \xi_{tu}), f(t)>0, t,u\in [0,1].$

We'll study the equation (\ref{e5}) in the class of
translational-invariant functions $f(t,x)$, i.e $f(t,x)=f(t)\in
C[0,1]$ for any $x\in V$ and it can be written as
\begin{equation}\label{e1.2}
(R_{k}f)(t)=f(t),
\end{equation}

Note that equation (\ref{e1.2}) is not linear for any $k\geq 1$.
For every $k\in\mathbb{N}$ we consider an integral operator
$H_{k}$ acting in the cone $C^{+}[0,1]$ i.e.,
$$(H_{k}f)(t)=\int^{1}_{0}K(t,u)f^{k}(u)du, \,\ k\in\mathbb{N}.$$

The operator $H_{k}$ is called Hammerstein's integral operator of
order $k$. Clearly, if $k\geq2$ then $H_{k}$ is a nonlinear
operator.

\begin{lemma}\label{l2.1.}\cite{ehr2013} Let $k\geq2$. The equation
\begin{equation}\label{e2.1}
R_{k}f=f, \,\ f\in C^{+}_{0}[0,1]
\end{equation}
has a nontrivial positive solution iff the Hammerstein's operator
has a positive eigenvalue, i.e. the Hammerstein's equation
\begin{equation}\label{e2.2} H_{k}f=\lambda f, \,\ f\in C^{+}[0,1]
\end{equation}
has a nonzero positive solution for some $\lambda>0$.
\end{lemma}

It is easy to check that if the number $\lambda_0>0$ is an
eigenvalue of the operator $H_{k}$, then an arbitrary positive
number is an eigenvalue of the operator $H_{k}$ (see Theorem 3.7
\cite{ehr2013}), where $k\geq2$. Consequently, we obtain

\begin{lemma}\label{l2.2.} Let $k\geq 2$. The equation (\ref{e2.1})
has a nontrivial positive solution iff the Hammerstein's operator
$H_{k}$ has a nontrivial positive fixed point, moreover
$N_{fix}^{+}(R_k)=N_{fix}^{+}(H_k)$, where $N_{fix}^{+}(T)$ is a
number of nontrivial positive fixed points of the operator $T$.
\end{lemma}

\section{Hammerstein's operator $H_{3}$ with degenerate kernel}

Let $\varphi_1(t), \,\ \varphi_2(t)$ and $\psi_1(t), \,\
\psi_2(t)$ are positive functions from $C_0^+[0,1]$. Suppose that
$\varphi_1(t)>0, \,\ \psi_1(t)>0.$ We consider Hammerstein's
operator $H_3:$

$$(H_3f)(t)=\int\limits\limits_0^1(\varphi_1(t)\psi_1(u)+\varphi_2(t)\psi_2(u))f^3(u)du$$
and cubic operator  $P$ on $\mathbb{R}^{2}$ by the rule

$$P(x, y)=(\alpha_{11}x^3+3\alpha_{12}x^2y+3\alpha_{21}xy^2+\alpha_{22}y^3, \,\  \beta_{11}x^3+3\beta_{12}x^2y+3\beta_{21}xy^2+\beta_{22}y^3).$$

Here

$$\alpha_{11} = \int\limits^{1}_{0}\psi_1(u)\varphi_1^3(u)du>0, \,\ \alpha_{12} = \int\limits^{1}_{0}\psi_1(u)\varphi_1^2(u)\varphi_2(u)du>0,$$
$$\alpha_{21} = \int\limits^{1}_{0}\psi_1(u)\varphi_1(u)\varphi_2^2(u)du>0, \,\ \alpha_{22} = \int\limits^{1}_{0}\psi_1(u)\varphi_2^3(u)du>0;$$

$$\beta_{11} = \int\limits^{1}_{0}\psi_2(u)\varphi_1^3(u)du>0, \,\ \beta_{12} = \int\limits^{1}_{0}\psi_2(u)\varphi_1^2(u)\varphi_2(u)du>0,$$
$$\beta_{21} = \int\limits^{1}_{0}\psi_2(u)\varphi_1(u)\varphi_2^2(u)du>0, \,\ \beta_{22} = \int\limits^{1}_{0}\psi_2(u)\varphi_2^3(u)du>0.$$

\begin{lemma}\label{l3.1.}The Hammerstein's operator
$H_3$ has a nontrivial positive fixed point iff the cubic
operator $P$ has a nontrivial positive fixed point, moreover
$N_{fix}^{+}(H_3)=N_{fix}^{+}(P)$.
\end{lemma}

\proof $(a)$ Put

$$\mathbb{R}_{2}^{+}=\{(x,y)\in \mathbb{R}^{2}: \,\ x\geq0, y\geq0\},$$
$$\mathbb{R}_{2}^{>}=\{(x,y)\in \mathbb{R}^{2}: \,\ x>0,
y>0\}.$$\\

Let the Hammerstein's operator $H_3$ has a nontrivial positive
fixed point $f(t)\in C_0^+[0,1]$. Let

\begin{equation}\label{e3.1}
c_1=\int\limits^{1}_{0}\psi_1(u)f^3(u)du
\end{equation}

and

\begin{equation}\label{e3.2}
c_2=\int\limits^{1}_{0}\psi_2(u)f^3(u)du.
\end{equation}

Clearly, $c_1>0, \,\ c_2>0$, i.e. $(c_1,
c_2)\in\mathbb{R}_{2}^{>}$. Then for the function $f(t)$ the
equality

\begin{equation}\label{e3.3}
 f(t)=c_1\varphi_1(t)+c_2\varphi_2(t)
\end{equation}
 holds.

Consequently, for parameters $c_1, c_2$ from the equality
(\ref{e3.1}) and (\ref{e3.2}) we have the two identity:

$$c_1=\alpha_{11}c_1^3+3\alpha_{12}c_1^2c_2+3\alpha_{21}c_1c_2^2+\alpha_{22}c_2^3,$$
$$c_2=\beta_{11}c_1^3+3\beta_{12}c_1^2c_2+3\beta_{21}c_1c_2^2+\beta_{22}c_2^3.$$

Therefore, the point $(c_1, c_2)$ is fixed point of the cubic
operator $P.$

$(b)$ Assume, that the point $(x_0, y_0)$ is a nontrivial positive
fixed point of the cubic operator $P,$ i.e. $(x_0, y_0)\in
\mathbb{R}_{2}^{+}\setminus\{\theta\}$ and numbers $x_0, y_0$
satisfies following equalities

$$\alpha_{11}x_0^3+3\alpha_{12}x_0^2y_0+3\alpha_{21}x_0y_0^2+\alpha_{22}y_0^3=x_0,$$
$$\beta_{11}x_0^3+3\beta_{12}x_0^2y_0+3\beta_{21}x_0y_0^2+\beta_{22}y_0^3=y_0.$$

Similarly, we can prove that the function
$f_0(t)=x_0\varphi_1(t)+y_0\varphi_2(t)$ is a fixed point of the
Hammerstein's operator $H_3$ and $f_0(t)\in C_0^+[0,1]$. This
completes the proof.
\endproof

\section{Positive fixed points of cubic operators in cone $\mathbb{R}_2^{+}$}

We define cubic operator (CO) $\mathcal{C}$
in cone of the space $\mathbb{R}^{2}$ by rule
$$ \mathcal{C}(x,y) = (a_{11}x^{3}+3a_{12}x^2y+3a_{21}xy^2+a_{22}y^{3}, \,\ b_{11}x^{3}+3b_{12}x^2y+3b_{21}xy^2+b_{22}y^{3}).$$

Clearly, an arbitrary nontrivial positive fixed points of the (CO) $\mathcal{C}$ is strictly positive. We denote by $N_{fix}^{>}(V)$
which the number of fixed points (CO) $\mathcal{C}$ belong to
$\mathbb{R}_{2}^{>}$ (belong to
$\mathbb{R}_{2}^{+}\setminus\{\theta\}$).
\begin{lemma}\label{l4.1.}\cite{MSSX} i) If the point $\omega=(x_0,y_0)\in \mathbb{R}_{2}^{+}$ is a fixed point of (CO) $\mathcal{C}$,
then $\omega \in\mathbb{R}_{2}^{>}$ and $\xi_0=\frac{y_0}{x_0}$ is
a root of the  algebraic equation
\begin{equation}\label{e4.1}
a_{22}\xi^{4}+(3a_{21}-b_{22})\xi^{3}+(3a_{12}-3b_{21})\xi^{2}+(a_{11}-3b_{12})\xi-b_{11}=0.
\end{equation}

ii) If the positive number $\xi_0$
is a root of the algebraic equation (\ref{e4.1}), then the
point $\omega_0=(x_0,\xi_{0}x_0)\in \mathbb{R}_{2}^{>}$ is fixed
point of (CO) $\mathcal{C}$, where
\begin{equation}\label{e4.2}
x_{0}=\frac{1}{(a_{11}+3a_{12}\xi_0+3a_{21}\xi_0^2+a_{22}\xi_0^3)^{1/2}}.
\end{equation}
\end{lemma}

We put

$$\mu_0=a_{22}, \,\ \mu_1=3a_{21}-b_{22}, \,\  \mu_2=a_{12}-b_{21}, \,\ \mu_3=a_{11}-3b_{12}, \,\ \mu_4=-b_{11}$$

and we define the  polynomial $P_4(\xi)$ of order four by

\begin{equation}\label{add} P_4(\xi)=\mu_0\xi^{4}+\mu_1\xi^{3}+3\mu_2\xi^2+\mu_3\xi+\mu_4. \end{equation}
The following result follows from Lemma \ref{e4.1}.
\begin{pro}\label{l4.1.}
 The number of positive fixed points of the operator (CO) $\mathcal{C}$ is equal to the number of positive roots of the polynomial  $P_4(\xi)$.
\end{pro}
The proof of this result follows, if the point $\xi_0$ is a positive root of the polynomial $P_4(\xi)$, then we can find the positive fixed point of the operator (CO) $\mathcal{C}$ with this formula
$$\omega_0=\left(\frac{1}{(a_{11}+3a_{12}\xi_0+3a_{21}\xi_0^2+a_{22}\xi_0^3)^{1/2}}, \frac{\xi_0}{(a_{11}+3a_{12}\xi_0+3a_{21}\xi_0^2+a_{22}\xi_0^3)^{1/2}}\right)$$
from this, it follows that if there are many positive roots of the polynomial $P_4(\xi)$, then there are as many positive fixed points of the operator (CO) $\mathcal{C}$.
\begin{lemma}\label{l4.4.}\cite{MSSX} The cubic operator (CO) $\mathcal{C}$ has no more than three fixed points, i.e. $1\leq N_{fix}^{>}(\mathcal{C})\leq3.$
\end{lemma}

\section{Non-uniqueness of positive fixed points of the Hammerstein's operator $H_3$}

We define two continuous positive functions $\zeta_{1}(t)$ and $\zeta_{2}(t)$ on $[0,1]$\\

$$\zeta_{1}(u)=\left\{\begin{array}{ll}
\frac{1}{2}+\sin2\pi u, \ \ \ \mbox{if} \ \ u\in[0,\frac{1}{2}]\\[2mm]
\frac{1}{2}, \ \ \ \mbox{if} \ \ u\in[\frac{1}{2},1] \\
\end{array}\right., $$
and
$$\zeta_{2}(u)=\left\{\begin{array}{ll}
\frac{1}{2}, \ \ \ \mbox{if} \ \ u\in[0,\frac{1}{2}]\\[2mm]
\frac{1}{2}-\sin2\pi u, \ \ \ \mbox{if} \ \ u\in[\frac{1}{2},1].\\
\end{array}\right. $$\\

For each positive numbers $a, b$ we define continuous positive
functions $F_{1}(t;a,b)$ and $F_2(t;a,b)$ on $[0,1]$

$$F_{1}(t;a,b)=\left\{\begin{array}{ll}
a\cos\pi t+b, \ \ \ \mbox{if} \ \ t\in[0,\frac{1}{2}]\\[2mm]
b, \ \ \ \mbox{if} \ \ t\in[\frac{1}{2},1] \\
\end{array}\right., $$\vskip 0.2truecm

$$ F_2(t;a,b)=\left\{\begin{array}{ll}
b, \ \ \ \mbox{if} \ \ t\in[0,\frac{1}{2}]\\[2mm]
-a\cos\pi t+b, \ \ \ \mbox{if} \ \ t\in[\frac{1}{2},1].\\
\end{array}\right.$$\vskip 0.2truecm

 For each positive numbers $a, b$ we denote

$$\tilde{K}(t,u;a,b)=\zeta_{1}(u)F_{1}(t;a,b)+\zeta_{2}(u)F_{2}(t;a,b), \,\ t,u\in
[0,1].$$

\begin{thm}\label{addp2} Let  $\tilde{K}(t,u;a,b)$ is kernel, then:\\

$(i)$ $ a<\frac{35(44+15\pi)}{318}b$ there exists a unique positive fixed point of the operator $H_{3}$ and is:
$$f(t)=\frac{1}{\sqrt{\frac{177a}{35\pi}+\frac{44+15\pi}{6\pi}b}}\zeta_{1}(t)+\frac{1}{\sqrt{\frac{177a}{35\pi}+\frac{44+15\pi}{6\pi}b}}\zeta_{2}(t).$$
\end{thm}

$(ii)$ $ a=\frac{35(44+15\pi)}{318}b$ there exist exactly two positive fixed points of the operator $H_{3}$.

$(iii)$ $ a>\frac{35(44+15\pi)}{318}b$ there exist exactly three positive fixed points of the operator $H_{3}$.

\begin{proof} $(i)$ At first we'll find coefficients of  $P_4(\xi)$ polynomial.
$$a_{11}=\int\limits_{0}^{1} F_{1}(u;a,b)\zeta_{1}^3(u)du=\int\limits_{0}^{\frac{1}{2}}(a\cos\pi u+b)\left(\frac{1}{2}+\sin2\pi u\right)^3du+$$
$$+\int\limits_{\frac{1}{2}}^{1}b\left(\frac{1}{2}\right)^{3}du=\frac{527}{280\pi}a+\frac{17}{12\pi}b+\frac{b}{2}\ . $$
Analogously we get
$$a_{12}=\int\limits_{0}^{1}F_{1}(u;a,b)
\zeta_{1}^2(u)\zeta_{2}(u)du=\frac{29a}{40\pi}+\frac{3b}{4\pi}+\frac{b}{4}, \ \
a_{21}=\int\limits_{0}^{1}F_{1}(u;a,b)
\zeta_{1}(u)\zeta^{2}_{2}(u)du=\frac{7a}{24\pi}+\frac{3b}{4\pi}+\frac{b}{4},$$
$$a_{22}=\int\limits_{0}^{1}F_{1}(u;a,b)
\zeta^{3}_{2}(u)du=\frac{a}{8\pi}+\frac{17b}{12\pi}+\frac{b}{2}.$$

 After short calculation we get $b_{11}=a_{22}, b_{12}=a_{21}, b_{21}=a_{12}, b_{22}=a_{11}$. Consequently, we have

$$P_4(\xi)=a_{22}\xi^{4}+(3a_{21}-a_{11})\xi^{3}+(a_{11}-3a_{12})\xi-a_{22}.$$
If we find the roots of this polynomial  $P_4(\xi)$,
$$a_{22}\xi^{4}+(3a_{21}-a_{11})\xi^{3}+(a_{11}-3a_{12})\xi-a_{22}=0$$
$$(\xi-1)(\xi+1)(a_{22}\xi^2+(3a_{21}-a_{11})\xi+a_{22})=0$$
$$\xi_1=-1, \xi_2=1, \,\  a_{22}\xi^2+(3a_{21}-a_{11})\xi+a_{22}=0$$

The quadratic equation $a_{22}\xi^2+(3a_{21}-a_{11})\xi+a_{22}=0$ has no roots when $D=(3a_{21}-a_{11})^2-4a_{22}^2<0$, ie $ a<\frac{35(44+15\pi)}{318}b$. So  $P_4(\xi)$ polynomial has unique positive root, the root is $\xi =1$. From Lemma (\ref{e4.2}) and (\ref{e3.3}), the fixed point of the Hammerstein operator $H_{3}$ is as follows
$$f(t)=\frac{1}{\sqrt{\frac{177a}{35\pi}+\frac{44+15\pi}{6\pi}b}}\zeta_{1}(t)+\frac{1}{\sqrt{\frac{177a}{35\pi}+\frac{44+15\pi}{6\pi}b}}\zeta_{2}(t).$$

 $(ii)$ We saw above that the roots of the polynomial $P_4(\xi)$ are the roots $\xi_1=-1, \xi_2=1$ and $a_{22}\xi^2+(3a_{21}-a_{ 11})\xi+a_{22}=0$ is equal to the roots of the quadratic equation. If we look at $D$, the discriminant of this quadratic equation
$$D=(3a_{21}-a_{11})^2-4a_{22}^2=(3a_{21}-a_{11}-2a_{22})(3a_{21}-a_{11}+a_{22})$$
If we replace the coefficients in the expression $3a_{21}-a_{11}-2a_{22}$, we will see that this expression is negative. $3a_{21}-a_{11}+a_{22}=0$ i.e. $a=\frac{35(44+15\pi)}{318}b$, then $D=0$. It follows that the quadratic equation  $a_{22}\xi^2+(3a_{21}-a_{11})\xi+a_{22}=0$ has one root. This root is positive, becouse $3a_{21}-a_{11}=-a_{22}$.  $a_{22}$ is positive, so $3a_{21}-a_{11}$ is negative. Then the operator $H_{3}$ has two positive fixed points.

$(iii)$This proof as above, i.e.  the roots of the polynomial $P_4(\xi)$ are the roots $\xi_1=-1, \xi_2=1$ and $a_{22}\xi^2+(3a_{21}-a_{ 11})\xi+a_{22}=0$ is equal to the roots of the quadratic equation. If $D>0$ then $(3a_{21}-a_{11}-2a_{22})(3a_{21}-a_{11}+a_{22})>0$. Obviously $(3a_{21}-a_{11}-2a_{22})$ is negative. So $(3a_{21}-a_{11}+a_{22})$ must be negative, i.e. $(3a_{21}-a_{11}+a_{22})<0$. If we replace the coefficients in the expression $(3a_{21}-a_{11}-2a_{22})$, we will see  $ a>\frac{35(44+15\pi)}{318}b$. Now we consider the expression $3a_{21}-a_{11}$. Since $(3a_{21}-a_{11}+a_{22})<0$ and $a_{22}>0$,  $3a_{21}-a_{11}<0$. Since $3a_{21}-a_{11}<0$ and $D>0$, it follows that roots of the quadratic equation  $a_{22}\xi^2+(3a_{21}-a_{11})\xi+a_{22}=0$ are positive. So
the operator $H_{3}$ has   three positive fixed points.
\end{proof}

\section*{Acknowledgements}
The work supported by the fundamental project (number: F-FA-2021-425)  of The Ministry of Innovative Development of the Republic of Uzbekistan.

\section*{Statements and Declarations}

{\bf	Conflict of interest statement:}
On behalf of all authors, the corresponding author states that there is no conflict of interest.

\section*{Data availability statements}
The datasets generated during and/or analysed during the current study are available from the corresponding author on reasonable request.

\end{document}